


\magnification\magstephalf

\hoffset=.5truein \advance\hoffset by -5.5truemm
\hsize=5.7truein
\voffset=.4375truein
\vsize=8.25truein

\parskip=6pt plus2pt minus2pt
\smallskipamount=6pt plus2pt minus2pt
\medskipamount=12pt plus4pt minus4pt
\bigskipamount=18pt plus6pt minus6pt
\abovedisplayskip=12pt plus4pt minus4pt
\belowdisplayskip=12pt plus4pt minus4pt

\normalbaselineskip=14pt plus1pt minus1pt
\normallineskiplimit=1pt
\normallineskip=1pt
\normalbaselines


\font\titlefont    = cmbx12 at 16pt
\font\twelverm     = cmr12
\font\ninerm       = cmr9
\font\sixrm        = cmr6
\font\nineit       = cmti9
\font\twelvebf     = cmbx12
\font\ninebf       = cmbx9
\font\sixbf        = cmbx6

\font\twelvebfi    = cmmib10 at 12pt
\font\eightbfi     = cmmib8
\font\ninei        = cmmi9
\font\sixi         = cmmi6
\font\twelvebfsy   = cmbsy10 at 12pt
\font\ninesy       = cmsy9      
\font\sixsy        = cmsy6

\font\tenmsa       = msam10
\font\ninemsa      = msam9
\font\sevenmsa     = msam7
\font\sixmsa       = msam6
\font\tenmsb       = msbm10
\font\ninemsb      = msbm9
\font\sevenmsb     = msbm7
\font\sixmsb       = msbm6
\font\teneuf       = eufm10   
\font\nineeuf      = eufm9
\font\seveneuf     = eufm7
\font\sixeuf       = eufm6
\font\teneus       = eusm10   
\font\twelveeusb   = eusb10 at 12pt 
\font\nineeus      = eusm9
\font\seveneus     = eusm7
\font\sixeus       = eusm6


\newfam\msafam
\textfont\msafam           =\tenmsa 
\scriptfont\msafam         =\sevenmsa

\newfam\msbfam
\def\bbb{\fam\msbfam}
\textfont\msbfam           =\tenmsb 
\scriptfont\msbfam         =\sevenmsb

\newfam\euffam
\def\euf{\fam\euffam}
\textfont\euffam           =\teneuf 
\scriptfont\euffam         =\seveneuf 

\newfam\eusfam
\def\eus{\fam\eusfam}
\textfont\eusfam           =\teneus 
\scriptfont\eusfam         =\seveneus

\def\ninepoint{
  \textfont0                 = \ninerm
  \scriptfont0               = \sixrm
  \textfont1                 = \ninei
  \scriptfont1               = \sixi
  \textfont2                 = \ninesy 
  \scriptfont2               = \sixsy
  \textfont\bffam            = \ninebf   
  \scriptfont\bffam          = \sixbf   
  \textfont\msafam           = \ninemsa
  \scriptfont\msafam         = \sixmsa
  \textfont\msbfam           = \ninemsb
  \scriptfont\msbfam         = \sixmsb
  \textfont\euffam           = \nineeuf
  \scriptfont\euffam         = \sixeuf
  \textfont\eusfam           = \nineeus
  \scriptfont\eusfam         = \sixeus
  \def\rm{\ninerm\fam0}%
  \def\it{\nineit\fam1}%
  \def\bf{\ninebf\fam\bffam}%
  \rm
  \normalbaselineskip=12pt%
  \normallineskiplimit=1pt%
  \normallineskip=1pt%
  \normalbaselines}

\def\boldtwelvepoint{
  \textfont0                 = \twelvebf
  \textfont1                 = \twelvebfi
  \scriptfont1                = \eightbfi
  \textfont2                 = \twelvebfsy
  \textfont\eusfam           = \twelveeusb
  \def\rm{\twelvebf\fam\bffam}%
  \normalbaselineskip=18pt%
  \normallineskiplimit=0pt%
  \normallineskip=1.2pt%
  \normalbaselines\rm}


\headline={\ifnum\pageno>1\ifodd\pageno\hfil{\ninerm\folio}%
  \else{\ninerm\folio}\hfil\fi\else\hfil\fi}
\footline={\hfil}
\def\makeheadline{\vbox to0pt{\vskip-34pt%
   \line{\vbox to8.5pt{}\the\headline}\vss}\nointerlineskip}

\def\titlelines#1\par{\leavevmode\bigskip\vskip\parskip\nointerlineskip
  \centerline{\vbox{\baselineskip=24pt\titlefont
  \halign{\hfil##\hfil\cr#1\crcr}}}\par}
\def\authorlines#1\par{\medskip\vskip\parskip\nointerlineskip
  \centerline{\vbox{\twelverm\baselineskip=18pt
  \halign{\hfil##\hfil\cr#1\crcr}}}\par}

\newif\ifheadingfeed
\def\heading#1\par{%
  \ifheadingfeed\vskip0pt plus.3\vsize\penalty-250\vskip0pt plus-.3\vsize
    \else\headingfeedtrue\fi
  \bigskip\vskip\parskip
  \centerline{\vbox{\boldtwelvepoint
  \halign{\hfil##\hfil\cr#1\crcr}}}%
  \smallskip}
\def\subheading#1\par{\removelastskip\smallskip
  \penalty-75\noindent{\boldtenpoint#1}\enspace}
\def\subsubheading#1\par{\penalty-50{\it#1}\enspace}
\long\def\proclaim#1\par#2\par{\ifdim\lastskip<\smallskipamount
  \removelastskip\smallskip\fi\penalty-100{\bf#1}\enspace
  {\it#2}\smallskip\penalty100}
\def\demo#1\par{\ifdim\lastskip<\smallskipamount\removelastskip\smallskip
  \penalty-50\fi{\it#1}\enspace}
\def\enddemo{\penalty200\hbox{}\nobreak\hfill\space\hbox{$\square$}\smallskip
  \penalty-75}
\def\item#1\par#2\par{{\advance\leftskip by 2\parindent
  \noindent\kern-\parindent\rlap{#1}\kern\parindent#2\par}}
\def\eqalign#1{\null\,\vcenter{\openup\jot
  \ialign{&\strut\hfil$\displaystyle##$%
          &$\displaystyle{}##$\hfil\cr#1\crcr}}\,}

\def\references{\heading\tenbf References\par\smallskip
  \parskip=2pt plus.7pt minus.7pt\tolerance=1000\overfullrule=0pt
  \frenchspacing\ninepoint}
\def\ref[#1]#2\par {\noindent\kern-1.8cm\hbox to0.8cm{[#1]\hfil}#2%
  \leftskip= 1.8cm\par\leftskip= 0cm}
\def\endlines#1{\medskip\vskip\parskip\nointerlineskip
  \centerline{\vbox{\hrule width 1in}}
  \medskip\nointerlineskip\vskip\parskip
  \centerline{\ninepoint\vbox{\halign{\hfil##\hfil\cr#1\crcr}}}}
\def\date{\ifcase\month\or January\or February\or March\or April\or May%
  \or June\or July\or August\or September\or October\or November\or December%
  \fi\space\number\day,\space\number\year}

\newif\ifshowlabel\showlabelfalse
\def\eqno#1{\leqno{%
  \ifshowlabel\llap{\hbox to.5in{\eightrm\string#1\hss}}\fi
  \hskip\parindent#1}}
\catcode`\@=11
\def\eqalignno#1{\displ@y\tabskip=\centering
  \halign to\displaywidth{\hfil$\@lign\displaystyle{##}$\tabskip=0pt
    &$\@lign\displaystyle{{}##}$\hfil\tabskip=\centering
    &\kern-\displaywidth\rlap{\hskip\parindent$\@lign##$}%
    \tabskip=\displaywidth\crcr#1\crcr}}
\catcode`\@=12

\hyphenation{mani-fold mani-folds}


\def\Ad{\mathop{\rm Ad}}
\def\ad{\mathop{\rm ad}}
\def\cont{{\rm cont}}

\def\Hom{\mathop{\rm Hom}}
\def\hom{\mathop{\rm hom}}
\def\ss{{\rm ss}}
\def\Tr{\mathop{\rm Tr}}

\def\U{{\rm U}}
\def\Spin{{\rm Spin}}
\def\SO{{\rm SO}}
\def\O{{\rm O}}
\def\SU{{\rm SU}}

\def\IC{{\bbb C}}
\def\IR{{\bbb R}}
\def\IZ{{\bbb Z}}

\def\eufg{{\euf g}}
\def\eufh{{\euf h}}
\def\eufk{{\euf k}}

\def\eufz{{\euf z}}
\def\eufhom{\mathop{\euf hom}}
\def\eufX{{\euf X}}

\def\eusA{{\eus A}}
\def\eusB{{\eus B}}
\def\eusE{{\eus E}}
\def\eusG{{\eus G}}
\def\eusGhat{\widehat{\eus G}}
\def\eusL{{\eus L}}
\def\eusO{{\eus O}}
\def\eusU{{\eus U}}
\def\eusV{{\eus V}}

\def\epsilon{\varepsilon}
\def\phi{\varphi}

\mathchardef\square="0803

\def\StabA#1{\Gamma_{#1}}
\def\StabB#1{\widehat\Gamma_{#1}}
\def\StabC#1{\Gamma_{#1}}
\def\nablaA{\nabla_{\!\!A}}
\def\mapright#1{\buildrel#1\over\longrightarrow}


\titlelines A Compactness Theorem for\cr Invariant Connections\cr

\authorlines Johan R\aa de

\bigskip

{\ninepoint\narrower\narrower

Necessary and sufficient conditions are given
for the Palais-Smale Condition~C to hold for the Yang-Mills
functional for invariant connections on a principal bundle
over a compact $n$-dimensional manifold.
The connections are assumed to be invariant
under the action of Lie group on the base manifold
such that all orbits of the action 
have dimension greater than or equal to $n{-}3$.

}


\heading \S1. The compactness theorem

It is well-known that the Palais-Smale Condition C
does not hold for the Yang-Mills functional on a principal bundle 
over a compact four-manifold.
According to the Uhlenbeck weak compactness theorem,
a Palais-Smale sequence will in general only have a subsequence
that converges on the complement of a finite set of points.
Moreover, even on the complement of these points,
the convergence is not as good as one would desire.
One only gets convergence with one derivative in $L^2$.
It is inconvenient to work with connections with one derivative in $L^2$
as there is no slice theorem for such connections.

It has been conjectured
that Condition C holds if the Yang-Mills functional is 
restricted to connections that are invariant under a group action
on the four-manifold, 
provided that all orbits of the action have dimension greater than or
equal to one.
The conjecture has been verified in a number of special cases;
see Examples 2 and 3 in \S2.
In this paper we settle this conjecture 
for manifolds of any dimension and group actions all of whose orbits have
codimension at most three.

Let $X$ be a smooth compact $n$-dimensional Riemannian manifold,
let $G$ be a compact Lie group and let $P$ be a smooth principal
$G$-bundle over $X$.
Let $\eusA$ be the space of smooth connections on $P$
and let $\eusG$ be the group of smooth gauge transformations 
(bundle automorphisms) of $P$.
Let 
$$\eusB=\eusA/\eusG$$
be the space of gauge equivalence classes of smooth connections on $P$.
We denote the gauge equivalence class of $A\in\eusA$ by $[A]\in\eusB$.

We choose a $G$-invariant positive definite 
inner product on the Lie algebra $\eufg$;
if $G$ is semisimple the negative of the Killing form provides
a canonical invariant inner product.
Then the Yang-Mills functional is given by ${\eus YM}[A]={1\over2}\int_X|F_A|^2dx$
where $F_A$ is the curvature of $A$.
The critical points of the Yang-Mills functional
are given by Yang-Mills equation $d_A^*F_A=0$.

Next we let $H$ be a compact Lie group
with a smooth action $\rho$ on $X$ by isometries
and with a lifted smooth action $\sigma$ on $P$ by bundle maps;
$$\matrix{
    H\times P & \mapright{\sigma} & P \cr
    \downarrow && \downarrow \cr
    H\times X & \mapright{\rho} & X. \cr
         }
$$
The action of $H$ on $P$ induces actions of $H$ on $\eusA$ and $\eusG$.
Let $\eusA^H$ and $\eusG^H$ be the fixed point sets of these actions.
These are the space of $H$-invariant connections 
and the group of $H$-invariant gauge transformations.
Let 
$$\eusE=\eusA^H/\eusG^H.$$
We denote the equivalence class of $A\in\eusA^H$ by $[A,\sigma]\in\eusE$.

There is a natural map 
$$\pi:\eusE\to\eusB$$
given by $\pi[A,\sigma]=[A]$.
Note that this map is not in general injective.
This occurs when $H$-invariant connections
are gauge equivalent, but not gauge equivalent through $H$-invariant
gauge transformation; see Example 1 in \S2.

A sequence $A_n\in\eusA$ is said to be a Palais-Smale sequence
if there exists $M>0$ such that
$$\displaylines{
    \|F_{A_n}\|_{L^2(X,\Lambda^2T^*X\otimes\Ad P)}  \le M
      \rlap{\qquad for all $n$} \cr
    \noalign{\line{and\hfill}}
    \big\|d^*_{A_n}F_{A_n}\big\|_{L^{-1,2}_{A_n}(X,T^*X\otimes\Ad P)}
       \to 0 \rlap{\qquad as $n\to \infty$.} \cr
          }
$$
Here
$\|\cdot\|_{L^{-1,2}_A}$ is the norm on $L^{-1,2}$
dual to the norm
$\|b\|_{L^{1,2}_A}^2
  = \|b\|_{L^2}^2
    + \|\nablaA b\|_{L^2}^2 
$
on $L^{1,2}$.
One says that Condition C holds if every Palais-Smale sequence 
has a subsequence that, after gauge transformations,
converges to a Yang-Mills connection.

As we have mentioned, invariant connections can be gauge equivalent
via non-invariant gauge transformations.
Therefore there are two forms of Condition C for the Yang-Mills
functional for invariant connections:
a weak form where we allow any gauge transformations
and a strong form where we require the gauge transformations to
be invariant.
The weak form suffices for minimizing the Yang-Mills functional over
the invariant connections;
the strong form is needed for Morse theory 
on the space of invariant connections.
Theorem 1 says that the weak form of Condition C always holds.
Theorem 2 gives a necessary and sufficient condition
for the strong form of Condition C to hold.
Theorem 3 provides practical means of verifying this condition 
in many cases of interest.

\proclaim Theorem 1.

If all orbits of the action of $H$ on $X$ have dimension $\ge n-3$
and if $A_n\in\eusA^H$ is a Palais-Smale sequence,
then there exists a subsequence, which we also denote $A_n$,
and a sequence $g_n\in\eusG$ such that
$g_n.A_n$ converge strongly in $L^{1,2}_3$ to a Yang-Mills
connection $A_\infty\in\eusA^H$ as $n\to\infty$.

\proclaim Theorem 2.

{
If all orbits of the action of $H$ on $X$ have dimension $\ge n-3$,
then the following two conditions are equivalent:

\item\rm(1)

For any Palais-Smale sequence $A_n\in\eusA^H$
there exists a subsequence, which we also denote $A_n$,
and a sequence $g_n\in\eusG^H$ such that
$g_n.A_n$ converge strongly in $L^{1,2}_3$ to a Yang-Mills
connection $A_\infty\in\eusA^H$ as $n\to\infty$.

\item\rm(2)

Any Yang-Mills connection in $\eusB$ has only finitely many preimages
in $\eusE$ under the map $\pi$.

}

\proclaim Theorem 3.

{
{\bf(a)}
If $G$ is semisimple, then any irreducible connection in $\eusB$ 
has only finitely many preimages in $\eusE$ under the map $\pi$.
If furthermore $H$ is connected or $Z(G)=1$, 
then there exist at most one preimage.

{\bf(b)}
If $H$ is semisimple, 
then any connection in $\eusB$ has only
finitely many preimages in $\eusE$ under the map $\pi$.

}

That a connection is irreducible here means that
its isotropy subgroup is isomorphic to the center $Z(G)$ of $G$.
The Morrey norm $L^{1,2}_3$ is discussed in \S3.
In two and three dimensions $L^{1,2}_3=L^{1,2}$.
In four and more dimensions connections in $L^{1,2}_3$ 
have as good analytical properties as 
connections in $L^{1,2}$ in three dimensions.

Theorem 1 and the implication (2)$\Rightarrow$(1) in Theorem 2 are
proven in \S4.
The implication (1)$\Rightarrow$(2) in Theorem 2
is proven in \S6.
Theorem 3 is a direct consequence of Proposition 6.2 and Lemma 5.5.
Proposition 6.2 gives a detailed description of the preimages under
$\pi$ of elements of $\eusB$.

T.~Parker, [21] Theorem 3.1,
has given a proof of Condition (1) in Theorem 2 under the
assumption that $X$ is four-dimensional,
the orbits have dimension greater than or equal to one,
$G$ is semisimple and $\eusA^H$ does not contain any reducible
connections.
However, U.~Gritsch [11] has pointed out a gap in Parker's proof;
in the proof of [21] Theorem 3.1 it is tacitly assumed,
with no justification,
that the gauge transformations given by [21] Theorem~2.3 are invariant.
Parker's philosophy is to work exclusively with invariant functions 
and sections.
Therefore it seems that to fill the gap one would need invariant good gauges.

{\it Acknowledgements.} The author whishes to thank Magnus Fontes,
Arne Meurmann, Jaak Peetre, Lorenzo Sadun and Karen Uhlenbeck 
for interesting discussions
and Gil Bor for comments on the manuscript.
He also whishes to thank
the Department of Mathematics at University of Texas at Austin,
where some of this research was carried out, and in particular
Dan Freed for their hospitality.


\heading \S2. Examples

Example 1 illustrates how invariant connections can be gauge equivalent
through noninvariant gauge transformations.
In examples 2 and 3 the weak form of Condition C is used to
minimize the Yang-Mills functional on $\eusE$.
In examples 4--6 the strong form of Condition C is
used to do Morse theory for the Yang-Mills functional.
Example 6 shows that one can do Morse theory even if the strong form
of Condition C fails, 
provided that the extent of the failure can be controlled.

{\it Example 1.}
Let $\SO(2)$ act on the trivial $\U(1)$-bundle over $S^1$ 
by letting it act in the natural way on the base
and trivially on the fibers.
The gauge equivalence class of a connection is uniquely
determined by its holonomy.
Hence $\eusB=\U(1)$.
The invariant connections are given by
$\nabla+i\alpha\,d\theta$ with $\alpha\in\IR$.
Here $\nabla$ denotes the trivial connection.
Thus $\eusA^{\SO(2)}=\IR$.
The invariant gauge transformations are given by constants,
so $\eusG^{\SO(2)}=\U(1)$.
These act trivially on $\eusA$,
so $\eusE=\eusA^{\SO(2)}/\eusG^{\SO(2)}=\eusA^{\SO(2)}=\IR$.
The connection $\nabla+i\alpha\,d\theta$ has holonomy $e^{-2\pi i\alpha}$.
It follows that $\eusE$ is the universal covering space of $\eusB$
and $\pi$ is the universal covering map.
We also see that two invariant connections $\nabla+i\alpha\,d\theta$
and $\nabla+i\beta\,d\theta$ are gauge equivalent 
if and only if $\alpha-\beta$ is an integer.
They are then related by the gauge transformation 
$e^{2\pi i(\alpha-\beta)\theta}$
which is not invariant unless $\alpha=\beta$.

{\it Example 2.} ([2], [21] \S4, [23], [26])
The spin-2 representation $\SO(3)\times\IR^5\to\IR^5$
induces an action of $\SO(3)$ on $S^4$ with the standard metric.
The generic orbit of this action is a copy of $\SO(3)/(\IZ_2\oplus\IZ_2)$.
There are two exceptional orbits that are copies of $\SO(3)/\O(2)=\IR P^2$.
Let $P_d$ be the principal $\SU(2)$-bundle over $S^4$
with second Chern number $d$.
The $\SO(3)$-action on $S^4$ induces a $\Spin(3)$-action on $S^4$.
The lifts of this action to $P_d$
are classified by non-negative integral weights $(w_-,w_+)$
such that $d=(w_-^2{-}w_+^2)/8$
and either $w_+=w_-=0$ or $w_+\equiv w_-\equiv 1$ (mod 2).
The Yang-Mills functional achieves its minimum on each
corresponding space $\eusE$ of invariant connections.
This follows from Theorem 1,
but was originally shown in [2] and [23] by dimensional reduction.
If $w_-\ge3$ and $w_+\ge3$, then the resulting Yang-Mills connection
is non-(anti-)selfdual.
If $(w_-,w_+)\ne(0,0),(1,1)$, then the resulting Yang-Mills connection
is irreducible.
This gives irreducible non-(anti-)selfdual Yang-Mills connections
on $P_d$ for all $d\ne\pm1$.

\vskip-\parskip
This example gives the simplest known construction of non-(anti-)selfdual
$\SU(2)$ Yang-Mills connections over $S^4$.
The existence of such connections on $P_0$ is originally 
due to L.~Sibner, R.~Sibner and K.~Uhlenbeck [24]
using a different invariant minimization scheme.
Whether there exist non-(anti-)selfdual Yang-Mills connections
on $P_1$ and $P_{-1}$ remains an open problem.

{\it Example 3.} ([2] \S3 Remark 4, [26])
The standard representation $\SO(3)\times\IC^3\to\IC^3$
induces an isometric action of $\SO(3)$ on $\IC P^2$
with the Fubini-Study metric.
The generic orbit is a copy of $\SO(3)/\IZ_2$.
There are two exceptional orbits:
a copy of $\SO(3)/\SO(2)=S^2$ and a copy of $\SO(3)/\O(2)=\IR P^2$.
Let $P_d$ be the principal $\SU(2)$-bundle over $\IC P^2$
with second Chern number $d$.
The lifts of the induced $\Spin(3)$-action to $P_d$
are classified by non-negative integral weights $(w_-,w_+)$
such that $d=(w_-^2{-}w_+^2)/4$ and either $w_-=0$ and $w_+\equiv 0$ (mod 2) or
$w_-\equiv w_+\equiv1$ (mod 2).
As in Example 2, the Yang-Mills functional achieves its minimum on each
corresponding space $\eusE$ of invariant connections.
If $w_-\ge3$ and $w_+\ge3$, then the resulting Yang-Mills connection
is non-(anti-)selfdual.
If $w_-\equiv w_+\equiv1$ (mod 2) and $(w_-,w_+)\ne(1,1)$,
then the resulting Yang-Mills connection is irreducible.
This gives irreducible non-(anti-)selfdual Yang-Mills connections
on $P_d$ for all even $d\ne\pm2$.
This example is closely related to the previous one;
$\IC P^2$ is an $\SO(3)$-equivariant double cover of
$S^4$ branched along one of the exceptional orbits.

{\it Example 4.} ([12])
There is a family $M_{2g}$ of spin four-manifolds obtained as 
$\IZ_2$-quotients of the product of $S^2$ with a genus $2g$ surface.
The natural action of $\SO(3)$ on $S^2$
induces an action of $\SO(3)$ on $M_{2g}$.
The generic orbit of this action is a copy of $\SO(3)/\SO(2) \cong S^2$. 
There is a family of exceptional orbits parametrized by $S^1$ of orbit type 
$\SO(3)/\O(2) \cong \IR P^2$. 
The action has a natural lift to the each chiral spinor bundle $S_\pm(M_{2g})$.
There are no reducible connections on $S_\pm(M_{2g})$ so
$\eusE$ is a Hilbert manifold.
By a min-max procedure for the Yang-Mills functional over $\eusE$, 
using Theorems 2 and 3,
Gritsch shows that for generic $\Spin(3)$-invariant metrics on 
$M_{2g}$ with $g \ne 1$ there exist irreducible 
non-(anti-)selfdual Yang-Mills connections on $S_{\pm}(M_{2g})$. 

{\it Example 5.} ([13], [26], [27])
Let $\SO(3)$ act diagonally on $S^2\times S^2$
with a metric which is given by round metrics on each $S^2$,
possibly with different radii.
The generic orbit of this action is a copy of $\SO(3)$.
There are two exceptional orbits that are copies of $S^2$;
these are given by pairs of identical points and pairs
of antipodal points on $S^2$.
Let $P_d$ be the principal $\SU(2)$-bundle over $S^2\times S^2$
with second Chern number $d$.
The lifts of the induced $\Spin(3)$-action are classified
by non-negative integral weights $(w_-,w_+)$ such that
$w_+\equiv w_-$ (mod 2) and $d=(w_-^2-w_+^2)/2$.
There are reducible connections in $\eusE$
so this space is not a Hilbert manifold.
Thus one introduces the Hilbert manifold $\widetilde\eusE$
of $\Spin(3)$-invariant connections with framings along an orbit.
Let $k$ be a nonnegative integer.
In [13] we use equivariant Morse theory 
for a perturbed Yang-Mills functional on
$\widetilde\eusE$ with $(w_-,w_+)=(2k,0)$ or $(0,2k)$ 
to construct irreducible non-(anti-)selfdual Yang-Mills connections 
on $P_d$, with $d=\pm2k^2$, of arbitrarily high energy.
Such connections have also been obtained,
for a different set of values of $d$,
by H.-Y.~Wang [27] using gluing techniques.
The non-(anti-)selfdual Yang-Mills connections
constructed in [26] Example 9.4 however are reducible.

{\it Example 6.}
The Hopf fibration gives a free action of $\SO(2)$ on $S^3$.
There is a unique principal $\SU(2)$-bundle over $S^3$
and the action has a unique lift to this bundle:
the trivial bundle and the trivial lift.
It follows from Proposition 6.2 that
the trivial connection in $\eusB$ has infinitely many
preimages in $\eusE$ under the map $\pi$.
It then follows from Theorem 2 that
the Yang-Mills functional on the space $\eusE$ 
of invariant connections does not satisfy the strong
form of Condition C.
This is also shown in [15] by dimensional reduction
to a Yang-Mills-Higgs type functional over $S^2$.

\vskip-\parskip
However, any nontrivial Yang-Mills connection is irreducible,
and has only finitely many preimages in $\eusE$ by Theorem 3(a).
Arguing as in the proof of Theorem 2,
one can prove the following Lemma:
{\it The strong form of Condition C holds
for any Palais-Smale sequence that has energy bounded from
below by a positive number.}

\vskip-\parskip
Using this one can still establish a Morse theory
for a perturbed Yang-Mills functional on $\widetilde\eusE$,
the Hilbert manifold of invariant connections with invariant
framings along an orbit.
Construct a perturbed Yang-Mills functional
on $\widetilde\eusE$ that has isolated non-degenerate critical points
and is equal to the Yang-Mills functional except
near the degenerate critical points.
Using the fact that the Yang-Mills functional on $\eusB$ 
satisfies Condition C,
one can show that the critical values of the perturbed Yang-Mills functional
on $\widetilde\eusE$ form a discrete set.
One gets a decomposition
of $\widetilde\eusE$ as a CW-complex.
Each critical point of the perturbed Yang-Mills functional
contributes one cell.
The critical value zero is attained at the preimages in $\widetilde\eusE$ 
of the trivial connection.
Each of these contributes one 0-cell.
Thus there are infinitely many 0-cells.
The Hilbert manifold $\widetilde\eusE$ is connected, 
so there then have to be infinitely many 1-cells.
However, it follows from the above Lemma that
any positive critical value is attained at only finitely
many critical points and therefore contributes only finitely many cells.
We conclude that there exist $\SU(2)$
Yang-Mills connections over $S^3$ of arbitrarily high energy.


\heading \S3. Analysis in Morrey spaces

In Theorems 1 and 2 we establish convergence in the Morrey norm $L^{1,2}_3$.
There are two reasons for using this norm.
First, an $H$-invariant function in $L^{1,2}(X)$ 
is automatically in $L^{1,2}_3(X)$ if all orbits of the action of $H$ on
$X$ have dimension $\ge n-3$.
Second, 
the space $L^{1,2}_3(X)$ has as good analytical properties as the space 
$L^{1,2}$ on a three-manifold,
see Lemmas 3.2--3.5, in particular Lemma 3.3.
In this section we review the basic properties of Morrey spaces.

The Morrey space $L^p_\lambda(\IR^n)$,
with $p\in[1,\infty)$ and $\lambda\in\IR$,
is defined as the space of all $f\in L^p(\IR^n)$
such that
$$\sup_{\rho\in(0,1]}\sup_{x\in\IR^n} \rho^{\lambda-n}
               \| f\|_{L^p(B_\rho(x))} ^p 
  <\infty .
$$
It is a Banach space with norm
$$\|f\|_{L^p_\lambda(\IR^n)}^p
  = \|f\|_{L^p(\IR^n)}^p
          + \sup_{\rho\in(0,1]}\sup_{x\in\IR^n}
               \rho^{\lambda-n} \|f\|_{L^p(B_\rho(x))} ^p .
$$
The Morrey space $L^{k,p}_\lambda(\IR^n)$,
with $k$ a positive integer, $p\in[1,\infty)$ and $\lambda\in\IR$,
is defined as the space
of all $f\in L^{k,p}(\IR^n)$ such that
$\partial^\alpha f\in L^p_\lambda(\IR^n)$ for all $\alpha$
with $|\alpha|\le k$.
It is a Banach space with norm
$$\|f\|_{L^{k,p}_\lambda(\IR^n)}^p
  = \sum_{|\alpha|\le k} \|\partial^\alpha f\|_{L^p_\lambda(\IR^n)}^p .
  \eqno(3.1)
$$
The Morrey space $f\in L^{-k,p}_\lambda(\IR^n)$,
with $k$ a positive integer, $p\in[1,\infty)$ and $\lambda\in\IR$,
is defined as the space of all $f\in L^{-k,p}(\IR^n)$
for which there exist $g_\alpha\in L^p_\lambda (\IR^n)$, $|\alpha|\le k$,
such that
$$f=\sum_{|\alpha|\le k} \partial^\alpha g_\alpha.
  \eqno(3.2)
$$
It is a Banach space with norm
$$\|f\|_{L^{-k,p}_\lambda(\IR^n)}^p = \inf_{(g_\alpha)}\sum_{|\alpha|\le k}
  \|g_\alpha\|_{L^p_\lambda(\IR^n)}^p
  \eqno(3.3)
$$
where we take infimum over all collections $(g_\alpha)_{|\alpha|\le k}$
of functions in $L^p_\lambda(\IR^n)$ that satisfy (3.2).
Note that
$$L^{k,p}_\lambda(\IR^n) = \left\{\eqalign{
    & 0 && \hbox{for $\lambda<0$} \cr
    & L^{k,\infty}(\IR^n) && \hbox{for $\lambda=0$} \cr
    & L^{k,p}(\IR^n) \qquad && \hbox{for $\lambda\ge n$.} \cr
                                   }\right.
$$
Thus Morrey spaces are of interest only for $\lambda\in(0,n)$.
For a brief introduction to Morrey spaces, see [9] Section 1 in Chapter 3.

Let $\eusU$ be an open subset of $\IR^n$.
Then the Morrey space $L^{k,p}_\lambda(\eusU)$,
with $k$ an integer, $p\in[1,\infty)$ and $\lambda\in\IR$,
is defined as the space of all $f\in L^{k,p}(\eusU)$
such that there exists $g\in L^{k,p}_\lambda(\IR^n)$
with $f=g\big|_\eusU$.
For negative $k$ the restriction is to be understood in the sense
of distributions.
It is a Banach space with norm
$$\|f\|_{L^{k,p}_\lambda(\eusU)}
  = \inf_g  \|g\|_{L^{k,p}_\lambda(\IR^n)}
  \eqno(3.4)
$$
where we take infimum over all $g\in L^{k,p}_\lambda(\IR^n)$
such that $f=g\big|_\eusU$.

Let $X$ be an $n$-dimensional smooth compact Riemannian manifold.
Then the Morrey space $L^p_\lambda(X)$,
with $k$ a positive integer and $\lambda\in\IR$,
is defined as the space
of all $f\in L^p(X)$ such that the norm
$$\|f\|_{L^p_{\lambda}(X)}^p
  =\|f\|_{L^p(X)}^p
   + \sup_{\rho\in(0,\rho_0/2]}\sup_{x\in X}
   \rho^{\lambda-n} \|f\|_{L^p(B_\rho(x))}^p 
  \eqno(3.5)
$$
is finite;
here $\rho_0$ is the injectivity radius of $X$.
This norm is invariant under isometries of $X$.
Let $E$ be a Euclidean vector bundle over $X$.
Let $A$ be a smooth metric connection on $E$.
Then the Morrey space $L^{k,p}_\lambda(X,E)$,
with $k$ a nonnegative integer, $p\in[1,\infty)$ and $\lambda\in\IR$,
is defined as the space of all $s\in L^{k,p}(X,E)$
such that the norm
$$\|s\|_{L^{k,p}_{\lambda,A}(X,E)}^p
  =\sum_{j=0}^k \bigl\|(|\nablaA^js|)\bigr\|_{L^p_\lambda(X)}^p
  \eqno(3.6)
$$
is finite.
The Morrey space $L^{-k,p}_\lambda(X,E)$,
with $k$ a positive integer, $p\in[1,\infty)$ and $\lambda\in\IR$,
is defined as the space of all $s\in L^{-k,p}(X,E)$
such that the norm
$$\|s\|_{L^{-k,p}_{\lambda,A}(X,E)}^p
  = \inf_{(t_j)}\sum_{j=0}^k
       \bigl\|(|t_j|)\bigr\|_{L^p_\lambda(X)}^p
  \eqno(3.7)
$$
is finite;
here we take infimum over all $(t_0,\ldots,t_k)$ such that
$t_j\in L^p_\lambda\bigl(X,(T^*X)^{\otimes j}\otimes E\bigr)$
and
$s = \sum_{j=0}^k (\nablaA^*)^jt_j$.
Different smooth connections $A$ give equivalent norms.

We can get equivalent norms for $L^{k,p}_\lambda(X,E)$
as follows.
Choose an atlas of good local coordinate charts
$\Phi_\nu:\eusU_\nu\to \IR^n$
and good lifts $\Psi_\nu:E\big|_{\eusU_\nu}\to\IR^n\times\IR^N$.
Here good means that the maps can be extended smoothly to
neighborhoods of $\overline\eusU_\nu$.
Then
$$\|s\|^p =
  \sum_\nu \bigl\| (\Psi^{-1}_\nu)^* s
                \|_{L^{k,p}_\lambda(\Phi_\nu(\eusU_\nu),\IR^N)}^p 
  \eqno(3.8)
$$
is a norm for $L^{k,p}_\lambda(X,E)$,
with $k\in\IZ$, $p\in[1,\infty)$ and $\lambda\in\IR$;
the details are left to the reader.

The norms (3.5)--(3.7) have the advantage of being 
invariant under isometries of $X$.
The norm (3.8) can be used to reduce 
analysis in Morrey spaces on manifolds to analysis
in Morrey spaces on $\IR^n$.

The following Lemma shows why Morrey spaces are useful
in gauge theory for invariant connections,
and possibly in other invariant situations as well.

\proclaim Lemma 3.1.

{
Let $H$ be a compact Lie group that acts smoothly
on $X$, in such a way that all $H$-orbits have dimension $\ge n-d$,
and that acts smoothly on $E$ covering the action on $X$.
If $s\in L^{k,p}(X,E)$, with $k\in\IZ$ and $p\in[1,\infty)$,
is $H$-invariant, then $s\in L^{k,p}_d(X,E)$.
If $A$ is an $H$-invariant connection on $E$,
then
$$\|s\|_{L^{k,p}_{d,A}(X)} \le  c\|s\|_{L^{k,p}_A(X)}.$$
The constant $c$ does not depend on $A$.
}

\demo Proof.

First one shows that if $\eusU$ is an open subset of $\IR^n$,
$\eusV$ is an open subset of $\IR^n$ 
with compact closure  contained in $\eusU$,
and $f\in L^p(\eusU)$ only depends on $(x_1,\ldots,x_d)$,
then $f\in L^p_d(\eusV)$ and
$$\|f\|_{L^p_d(\eusV)} \le c \|f\|_{L^p(\eusU)}.$$
The idea is that any ball of radius $r$ in $\eusV$
has the order of magnitude $r^{d-n}$ disjoint translates in $\eusU$
in the $(x_{d+1},\ldots,x_n)$ directions
and $f$ has the same $L^p$-norm on all these balls;
we leave the detailed verification to the reader.

Next consider the case in the Lemma with $k=0$.
In a neighborhood of each point on $X$ we can find local coordinates
where points with the same $(x_1,\ldots,x_d)$ lie in the same $H$-orbit.
The Lemma then follows by applying the above estimate to these
coordinate charts.

For $k>0$ the Lemma follows by applying the case $k=0$
to the functions $|\nablaA^js|$.
For $k<0$ the Lemma follows by applying the case $k=0$
to the functions $|t_j|$,
once we observe that by averaging we may take the sections $t_j$
in (3.7) to be $H$-invariant.
\enddemo

Next we review embedding theorems and elliptic estimates in Morrey spaces.
First we have an analogue of H\"older's inequality:

\proclaim Lemma 3.2.

Multiplication gives bounded linear maps
$$L^p_\lambda(X,E)\times L^q_\lambda(X,F) \to L^r_\lambda(X,E\otimes F)$$
for all $1\le p,q,r<\infty$ such that
$1/p+1/q=1/r$.

Using the norms (3.8)
we see that it suffices to show that multiplication gives bounded linear maps
$L^p_\lambda(\IR^n)\times L^q_\lambda(\IR^n)\to L^r_\lambda(\IR^n)$.
That follows immediately from H\"older's inequality applied to $\IR^n$
and to the balls $B_\rho(x)$.

Next we have that the Morrey spaces $L^{k,p}_d$ in $n$ dimensions
have embedding properties analogous to those of
the Sobolev spaces $L^{k,p}$ in $d$ dimensions.

\proclaim Lemma 3.3.

Let $\lambda\in(0,n)$.
Then there are embeddings
$$L^{1,p}_\lambda(X,E)\to  \cases{
              L^{p^*}_\lambda(X,E) & for $p\in(1,\lambda)$ \cr
               C^{0,\alpha}(X,E) & for $p\in(\lambda,\infty)$ \cr
             }
$$
where $1/{p^*}=1/p-1/\lambda$ and $\alpha=1-\lambda/p$.

Again it suffices to establish the analogous
embedding on $\IR^n$.
A  weaker version of the first embedding,
$L^{1,p}_\lambda(\IR^n)\to L^q_\lambda(\IR^n)$ for $q\in[1,p^*)$,
which suffices for our purposes,
was first shown by S.~Campanato, [4] Theorem~1.2.
The sharp embedding $L^{1,p}_\lambda(\IR^n)\to L^{p^*}_\lambda(\IR^n)$
is due to D.R.~Adams [1] Theorem~3.2; see also [6]~Theorem~2.
The embedding
$L^{k,p}_\lambda(\IR^n)\to C^{0,\alpha}(\IR^n)$
is a classical result by C.B.~Morrey, [18] pp.~12--14;
see also [9] Theorem~1.2 in Chapter 3 and [10] Theorem 7.19.

From Lemma 3.3 we in particular get the embeddings
$$\eqalignno{
     L^{1,2}_3&\to L^6_3 \cr
    \noalign{\line{and\hfill}}
    L^{2,2}_3&\to L^{1,6}_3 \to C^{0,1/2} . \cr
               }
$$
By combining Lemmas 3.2 and 3.3 we get the multiplications
$$\eqalignno{
     & L^{2,2}_3 \times L^{2,2}_3 \to L^{2,2}_3 \cr
     & L^{2,2}_3 \times L^{1,2}_3 \to L^{1,2}_3  \cr
     \noalign{\line{and\hfill}}
     & L^{1,2}_3 \times L^{1,2}_3 \to L^3_3 \to L^2 . \cr
               }
$$
It follows that the group operations in $\eusG$,
the action of $\eusG$ on $\eusA$,
and the Yang-Mills functional
are continuous in the $L^{2,2}_3$-topology on $\eusG$ and
$L^{1,2}_3$-topology on $\eusA$.

There is also a compact embedding theorem.
This follows by standard methods from [1] Theorem 3.1.
For the convenience of the reader we have included a proof.

\proclaim Lemma 3.4.

Let $\lambda\in(0,n)$ and $p\in(1,\lambda)$.
Let $q\in[1,p^*)$ where $1/p^*=1/p-1/\lambda$.
Then the embedding $L^{1,p}_\lambda(X,E)\to  L^q_\lambda(X,E)$
is compact.

\demo Proof.

It suffices to show that if $\eusU$
is a bounded open subset of $\IR^n$,
then the embedding $L^{1,p}_\lambda(\eusU)\to  L^q_\lambda(\eusU)$
is compact.
We may assume that $q\in[p,p^*)$.
Let $\alpha=\lambda/q-\lambda/p^*\in(0,1]$.
We say that
$f\in C^{0,\alpha} L^q_\lambda(\IR^n)$
if $f\in L^q_\lambda(\IR^n)$ and
$$\sup_{0<|h|\le 1} |h|^{-\alpha}
         \bigl\| f(\cdot+h)-f(\cdot)\bigr\|_{L^q_\lambda(\IR^n)} < \infty.
$$
This is a Banach space with norm
$$\|f\|_{C^{0,\alpha} L^q_\lambda(\IR^n)}^q
  = \|f\|_{L^q_\lambda(\IR^n)}^q
    + \sup_{0<|h|\le 1} |h|^{-\alpha q}
            \bigl\| f(\cdot+h)-f(\cdot)\bigr\|_{L^q_\lambda(\IR^n)}^q .
$$
If $f\in L^{1,1}(\IR^n)$,
then
$$f(x)=c_n^{-1}\sum_{j=1}^n\int_{\IR^n}
       (x_j-y_j)\,|x-y|^{-n}\,\partial_j f(y)\,dy
$$
almost everywhere, where $c_n$ is the $(n-1)$-dimensional
measure of the unit sphere in $\IR^n$.
Thus
$$\eqalign{
  & |h|^{-\alpha}\bigl(f(x+h)-f(x)\bigr) \cr
  & \quad = c_n^{-1} \sum_{j=1}^n \int_{\IR^n} |h|^{-\alpha} \bigl(
      (x_j-y_j+h_j)\,|x-y+h|^{-n} - (x_j-y_j)\,|x-y|^{-n}
                                            \bigr) \partial_j f(y)\,dy \cr
          }
$$
almost everywhere. Now
$$|h|^{-\alpha}\bigl|(z_j+h_j)\,|z+h|^{-n} - z_j|z|^{-n}\bigr|
  \le c \bigl(|z+h|^{1-n-\alpha}+|z|^{1-n-\alpha}\bigr)
$$
almost everywhere. To see this, note that both sides are homogeneous with
respect to $(z,h)$ of degree $1{-}n{-}\alpha$.
Therefore it suffices to verify the inequality on the sphere $|z|^2+|h|^2=1$,
which is straightforward.
It follows that
$$\eqalign{
  |h|^{-\alpha}\bigl|f(x+h)-f(x)\bigr|
  & \le c\int_{\IR^n} \bigl(
      |x-y+h|^{1-n-\alpha} + |x-y|^{1-n-\alpha}
                          \bigr) |\,\nabla f(y)| \,dy \cr 
  & = c\int_{\IR^n}
      |x-y|^{1-n-\alpha} \bigl( |\nabla f(y)| + |\nabla f(y+h)| \bigr)\,dy \cr
          }
$$
almost everywhere. It now follows from [1] Theorem 3.1 or [6] Theorem 2
that the right hand side is bounded in $L^q_\lambda(\IR^n)$ uniformly in $h$.
We conclude that there is a continuous embedding
$L^{1,p}_\lambda(\IR^n)\to C^{0,\alpha}L^q_\lambda(\IR^n)$.

We define $C^{0,\alpha} L^q_\lambda(\eusU)$
by extension to $\IR^n$ as in (3.4).
It then follows that there is a continuous  embedding
$L^{1,p}_\lambda(\eusU)\to C^{0,\alpha}L^q_\lambda(\eusU)$.
It is not hard to show that if $\eusU$ is bounded,
then the embedding
$C^{0,\alpha} L^q_\lambda(\eusU)\to L^q_\lambda(\eusU)$ is compact.
\enddemo

\proclaim Lemma 3.5.

Let $p\in(1,\infty)$ and $k\in\IZ$.
Then any elliptic partial differential operator
$L^{k+m,p}_\lambda(X,E)\to L^{k,p}_\lambda(X,F)$
of order $m$ with smooth coefficients is a Fredholm operator.

To prove this one has to establish interior estimates for
elliptic operators on bounded open subsets of $\IR^n$.
That was first done by S.~Campanato, [5] Theorem~10.1, in the case $p=m=2$,
which suffices for our purposes.
To do the general case,
one has to show that singular integral operators
are bounded $L^p_\lambda(\IR^n)\to L^p_\lambda(\IR^n)$.
That is a result due to J.~Peetre [22] Theorem~1.1;
see also [6] Theorem~3 and [25] Proposition~3.3.

{\it Remark 3.6.}
No one has bothered to show that the theory of
elliptic boundary value problems extends to Morrey spaces.
However, it is clear that the classical treatment of the Neumann and Dirichlet
problems for the Laplacian by even and odd reflection
across the boundary carries over to Morrey spaces.

{\it Remark 3.7.}
We have chosen to work with smooth connections
and smooth gauge transformations.
One could also complete the space of connections
and the group of gauge transformations in $L^{1,2}_3$
and $L^{2,2}_3$ respectively.
Then one should observe that
$C^\infty_0(\IR^n)$ is not dense in $L^{k,p}_\lambda(\IR^n)$
for $\lambda\in(0,n)$.
The closure of $C^\infty_0(\IR^n)$ in $L^{k,p}_\lambda(\IR^n)$
for $\lambda\in(0,n)$
is a closed subspace $C^0L^{k,p}_\lambda(\IR^n)$
of $L^{k,p}_\lambda(\IR^n)$.
It consists of all $f\in L^{k,p}_\lambda(\IR^n)$
such that $f(\cdot+h)\to f(\cdot)$ in $L^{k,p}_\lambda(\IR^n)$
as $h\to0\in\IR^n$.
It is of course a Banach space with the same norm as $L^{k,p}_\lambda(\IR^n)$.
By the usual techniques, one can define
$C^0 L^{k,p}_\lambda(X,E)$,
which then is the closure of $C^\infty(X,E)$ in $L^{k,p}_\lambda(X,E)$.
Thus the completions of $\eusA$ and $\eusG$ would be the space
of connections in $C^0L^{1,2}_3$
and the group of gauge transformations in $C^0L^{2,2}_3$.


\heading \S4. Invariant connections

When we stated the main theorems in \S2,
we considered a fixed lift $\sigma:H\times P\to P$
of the action $\rho:H\times X\to X$.
However, in the proofs of these theorems we will need to consider
other lifts.
In particular, the gauge transforms of the lift $\sigma$ will play
an important role.
The group $\eusG$ acts on the set of lifts of $\rho$ as follows:
$$(g.\sigma)(h)=g\,\sigma(h)\,g^{-1}.$$
Note that if $A$ is invariant under $\sigma$,
then $g.A$ is invariant under $g.\sigma$.
Also note that
$$g.\sigma=\sigma
  \quad\Leftrightarrow\quad
  \hbox{$g\circ\sigma(h)=\sigma(h)\circ g$ for all $h\in H$}
  \quad\Leftrightarrow\quad
  g\in\eusG^H,
  \eqno(4.1)
$$
i.e.~$\sigma$ is invariant under $g$
if and only if $g$ is invariant under $\sigma$.

If we wish to consider other lifts of $\rho$,
then it is natural to introduce the extended gauge group $\eusGhat$ defined
as the group of ordered pairs $(\phi,h)$
such that $h\in H$ and $\phi:P\to P$ is a bundle map
that covers the isometry $\rho(h):X\to X$.
There is a short exact sequence:
$$1\longrightarrow\eusG\longrightarrow\eusGhat
   \longrightarrow H\longrightarrow 1.
  \eqno(4.2)
$$
A lift $\sigma'$ of $\rho$ is by definition the the same as a splitting 
$H\to\widehat\eusG$ of this short exact sequence.

There are compatible group actions:
$$\eqalign{
    \eusG\times\eusA&\to\eusA\cr
    \eusGhat\times\eusA&\to\eusA\cr
    H\times\eusB&\to\eusB.\cr
          }
$$
The first two actions are defined the obvious way.
The third action was introduced in [8] and is defined as follows:
Given $h\in H$ and $[A]\in\eusB$,
we choose any bundle map $\phi$ that covers the action of $h$ on $X$.
Then we define $h.[A]=[\phi.A]\in\eusB$.
As $A$ and $\phi$ are unique up to gauge transformations,
$\phi.A$ is unique up to gauge transformations.
Hence $[\phi.A]$ is a well-defined element of $\eusB$.
We denote the fixed point set of the action $H\times\eusB\to\eusB$
by $\eusB^H$.

For any $A\in\eusA$ these actions have isotropy subgroups
$$\eqalign{
    \StabA{A} & = \bigl\{g\in\eusG \,|\, g.A=A\bigr\} \cr
    \StabB{A} & = \bigl\{(\phi,h)\in\eusGhat \,|\, \phi.A=A\bigr\} \cr
    \StabC{[A]} & = \bigl\{h\in H  \,|\, h.[A]=[A]\bigr\} . \cr
          }
$$
These are compact Lie groups,
and they form a short exact sequence
$$1\longrightarrow\StabA{A}\longrightarrow\StabB{A}
   \longrightarrow \StabC{[A]}\longrightarrow 1.
  \eqno(4.3)
$$
We have $[A]\in\eusB^H$ if and only if $\StabC{[A]}=H$.

If $A$ is invariant under any lift $\sigma'$ of $\rho$,
then $[A]\in\eusB^H$.
Conversely, if $[A]\in\eusB^H$,
then the lifts under which $A$ is invariant are by definition precisely the 
splittings $H\to\StabB{A}$ of the short exact sequence
$$1\longrightarrow\StabA{A}\longrightarrow\StabB{A}
   \longrightarrow H\longrightarrow 1.
  \eqno(4.4)
$$

Arguing as usual,
but using the Morrey space estimates of \S3 instead of the usual Sobolev
space estimates,
we see that for $\epsilon>0$ small enough
$$\eusO_{A_0,\epsilon}
  = \bigl\{ A_0+a\in\eusA
    \,\big|\, \hbox{
          $d^*_{A_0}a=0$ and
          $\|a\|_{L^{1,2}_{3,A_0}(X,T^*X\otimes \Ad P)} <\epsilon$}
    \bigr\}
$$
is a local slice through $A_0$ for the action of $\eusG$ on $\eusA$.

\proclaim Lemma 4.1.

For any $A_0\in\eusA$ with $[A_0]\in\eusB^H$
there exists $\epsilon>0$ such that
$$\StabB{A}\subseteq \StabB{A_0}$$
for all $A\in\eusO_{A_0,\epsilon}$.

\demo Proof.

First I claim that
the action $\eusGhat\times\eusA\to\eusA$
restricts to an action
$$\StabB{A_0}\times\eusO_{A_0,\epsilon}\to\eusO_{A_0,\epsilon}.$$
In fact,
let $(\phi,h)\in\StabB{A_0}$ and $A=A_0+a\in\eusO_{A_0,\epsilon}$.
Then $\phi.A=A_0+\phi.a$.
The norm $\|\cdot\|_{L^{1,2}_{3,A_0}}$ is $\StabB{A_0}$-invariant,
so $\|\phi.a\|_{L^{1,2}_{3,A_0}}=\|a\|_{L^{1,2}_{3,A_0}}<\epsilon$.
Moreover,
$$d^*_{A_0}(\phi.a) =d^*_{\phi.A_0}(\phi.a)  = \phi.d^*_{A_0}a = 0 .
$$
The claim follows.

Let now $A\in\eusO_{A_0,\epsilon}$ and $(\phi,h)\in\StabB{A}$.
We have to show that $(\phi,h)\in\StabB{A_0}$.
That $A_0\in\eusB^H$ means that $\StabC{[A_0]}=H$.
It then follows from the exactness of the sequence (4.3) for $A_0$
that there exists a bundle map $\psi$ such that $(\psi,h)\in\StabB{A_0}$.
It follows from the exactness of (4.2) that $\psi=g\phi$
for some $g\in\eusG$.
Thus
$$(g\phi,h)\in\StabB{A_0}.
$$
Hence $g\phi.A\in\eusO_{A_0,\epsilon}$.
But $\phi.A=A$ so $g.A\in\eusO_{A_0,\epsilon}$.
It is well-known that,
for $\epsilon$ small enough,
$A,g.A\in\eusO_{A_0,\epsilon}$ implies $g\in\StabA{A_0}$.
Hence
$$(g,1)\in\StabB{A_0}.$$
We conclude that $(\phi,h)\in\StabB{A_0}$.
\enddemo

We can now show an important semicontinuity property:

\proclaim Proposition 4.2.

For any $A_0\in\eusA$ there exists $\epsilon>0$ 
such that if $A\in\eusO_{A_0,\epsilon}$, $\sigma'$~is a lift of $\rho$,
and $A$ is invariant under $\sigma'$,
then $A_0$ is invariant under $\sigma'$.

\demo Proof.

If $[A_0]\notin\eusB^H$, then 
there exists $\epsilon>0$ such that $\eusO_{A_0,\epsilon}$
is disjoint from $\eusB^H$,
for $\eusB^H$ is a closed subset of $\eusB$.
It then follows that no connection in $\eusO_{A_0,\epsilon}$ 
is invariant under any lift of $\rho$, 
and the Proposition holds trivially.

If $[A_0]\in\eusB^H$, then it follows from Lemma 4.1
that, for $\epsilon$ small enough, $\StabB{A}\subseteq\StabB{A_0}$.
A lift $\sigma'$ that preserves $A_0$ is by definition
the same as a continuous right inverse of the homomorphism $\StabB{A_0}\to H$.
A lift $\sigma'$ that preserves $A$ is by definition
a continuous right inverse of $\StabB{A}\to H$.
The Proposition follows.
\enddemo

The following Proposition is essentially due to M.~Furuta, [8] Lemma 2.2;
see also [3] Proposition 4.3.

\proclaim Proposition 4.3.

$\pi(\eusE)$ is a closed subset of $\eusB$.

\demo Proof.

Assume that $A_n\in\eusA^H$ and
$\pi[A_n,\sigma]=[A_n]\to[A_\infty]$ in $L^{1,2}_3$ as $n\to\infty$.
This means that there exist $g_n\in\eusG$ such that 
$g_n.A_n\to A_\infty$ in $L^{1,2}_3$ as $n\to\infty$.
We have to show that $[A_\infty]\in\pi(\eusE)$.

Let $\epsilon$ be as in Proposition 4.2.
We may assume that $g_n.A_n\in\eusO_{A_\infty,\epsilon}$ for large $n$.
As $A_n$ is invariant under $\sigma$,
we have that $g_n.A_n$ is invariant under $g_n.\sigma$.
It then follows from Proposition 4.2 that
$A_\infty$ is  invariant under $g_n.\sigma$.
Hence $g_n^{-1}.A_\infty$ is invariant under $\sigma$,
i.e.~$g_n^{-1}.A_\infty\in\eusA^H$ for large $n$.
Thus $[A_\infty]=[g_n^{-1}.A_\infty]=\pi[g_n^{-1}.A_\infty,\sigma]$.
\enddemo

\demo Proof of Theorem 1.

Let $A_n\in\eusA$ be a Palais-Smale sequence with
$[A_n]\in\eusB^H$.
Then $A_n$ is invariant under the action of $\StabB{A_n}$ on $P$.
It follows from Lemma 3.1 that
$$\|F_{A_n}\|_{L^2_3(X,\Lambda^2T^*X\otimes\Ad P)}
  \le c \|F_{A_n}\|_{L^2(X,\Lambda^2T^*X\otimes\Ad P)}
$$
and
$$\big\|d^*_{A_n}F_{A_n}\big\|_{L^{-1,2}_{3,A_n}(X,T^*X\otimes\Ad P)}
  \le c \big\|d^*_{A_n}F_{A_n}\big\|_{L^{-1,2}_{A_n}(X,T^*X\otimes\Ad P)} 
$$
uniformly in $n$.
Hence
$$\|F_{A_n}\|_{L^2_3(X,\Lambda^2T^*X\otimes\Ad P)}
  \le c M \rlap{\qquad for all $n$}
$$
and
$$\big\|d^*_{A_n}F_{A_n}\big\|_{L^{-1,2}_{3,A_n}(X,T^*X\otimes\Ad P)}
  \to 0 \rlap{\qquad as $n\to\infty$.}
$$
It follows
that there exists a sequence of gauge transformations 
$$g_n\in\eusG$$
such that $g_n.A_n$ converge in $L^{1,2}_3$-norm
to a smooth Yang-Mills connection $A_\infty$ as $n\to\infty$.
In fact, for $n=3$, where $L^{1,2}_3=L^{1,2}$, this is well-known.
It is proven exactly the same way for $n\ge4$,
using the Morrey space estimates in \S3 instead of the usual
Sobolev space estimates.

It follows from Proposition 4.3 that
$[A_\infty]\in\pi(\eusE)$.
In other words,
$A_\infty$ is gauge equivalent to a connection in $\eusA^H$.
After adjusting the gauge transformations $g_n$ we may thus assume that
$A_\infty\in\eusA^H$.
\enddemo

\demo Proof of Theorem 2, (2)$\Rightarrow$(1).

Assume that Condition (2) holds.
Let $A_n\in\eusA^H$ be a Palais-Smale sequence.
It follows from Theorem 1 that there exists a sequence $g'_n\in\eusG$
and a Yang-Mills connection $A'_\infty\in\eusA^H$
such that
$$g'_n.A_n\to A'_\infty \rlap{\qquad in $L^{1,2}_3$\quad as $n\to\infty$.}
  \eqno(4.5)
$$
Let $\epsilon$ be as in Proposition 4.2.
We can choose the sequence $g'_n$
so that $g'_n.A_n\in\eusO_{A'_\infty,\epsilon}$
for large $n$.
We have that $g'_n.A_n$ is invariant under $g'_n.\sigma$.
It then follows from Proposition 4.2 that $A'_\infty$ is invariant under
$g'_n.\sigma$ for large $n$.
Hence ${g'_n}^{-1}.A_n$ is invariant under $\sigma$;
in other words
$${g'_n}^{-1}.A'_\infty\in\eusA^H
  \eqno(4.6)
$$
for large $n$.
In particular $[{g'_n}^{-1}.A'_\infty,\sigma]\in\pi^{-1}[A'_\infty]$
for large $n$.
It follows from Condition (2) that $\pi^{-1}[A'_\infty]$ is finite.
After passing to a subsequence we may therefore assume that
all $[{g'_n}^{-1}.A'_\infty,\sigma]$ are equal.
This means that for any integer $n_0$ there exists a sequence
$$g_n\in\eusG^H$$
such that
$g_n {g'_n}^{-1}.A'_\infty={g'_{n_0}}^{\!\!\!-1}.A'_\infty$ for all $n$.
Then
$g'_{n_0}g_n {g'_n}^{-1}\in\StabA{A'_\infty}$.
The group $\StabA{A'_\infty}$ is compact.
Hence, after passing to a subsequence,
$$g'_{n_0}g_n  {g'_n}^{-1}\to s
   \rlap{\qquad in $L^{2,2}_3$\quad as $n\to\infty$}
  \eqno(4.7)
$$
for some
$$s\in\StabA{A'_\infty}.
  \eqno(4.8)
$$
It follows from (4.7), (4.5), (4.8) and (4.6) that
$$g_n.A_n = (g_n {g'_n}^{-1}).(g'_n.A_n)
  \;\to\;
  {g'_{n_0}}^{\!\!\!-1} s.A'_\infty
  = {g'_{n_0}}^{\!\!\!-1}A'_\infty\in\eusA^H
$$
in $L^{1,2}_3$ as $n\to\infty$.
Condition (1) follows with
$A_\infty={g'_{n_0}}^{\!\!\!-1}.A'_\infty$.
\enddemo


\heading \S5. Homomorphisms of compact Lie groups

In this section we review some facts about homomorphisms of compact
Lie groups that will be used in the proofs of the 
implication (1)$\Rightarrow$(2) in Theorem 2 and Theorem 3(b).
Note that we do not assume the groups to be connected.

If $H$ is a compact Lie group and $K$ is a Lie group,
then we let $\Hom(H,K)$ denote the set of continuous,
and hence smooth,
homomorphisms $H\to K$.
We will view $\Hom(H,K)$ as a subset of the Banach manifold $C(H,K)$
of continuous maps $H\to K$.
The group $K$ acts on $\Hom(H,K)$ and $C(H,K)$ by conjugation.

The following rigidity theorem,
for compact Lie groups $H$,
was announced and a proof was outlined
by A.~Nijenhuis and R.W.~Richardson [20] Theorem C.
The details of the proof were carried out,
for compact topological groups $H$,
by D.H.~Lee [17] Theorem 2.
The proof is based on an idea by A.~Weil [28].
For the convenience of the reader
we give a self-contained proof essentially following [20].
Although we only need the rigidity theorem for compact Lie groups $H$,
we state and prove it for compact topological groups $H$
as that does not require any extra work.
For more information on the continuous cohomology
introduced in the proof, see [7], [19], [16] and [14]~Chapter~3.

\proclaim Lemma 5.1.

If $H$ is a compact topological group and $K$ is a Lie group,
then $\Hom(H,K)$ is a closed discrete union of $K$-orbits in $C(H,K)$.

\demo Proof.

It is clear that $\Hom(H,K)$ is a closed subset of $C(H,K)$.
We will now show that each $K$-orbit in $\Hom(H,K)$ is isolated.
We write $1$ for the identity elements in $H$ and $K$.
We also write 1 for the map $H\times H\to K$
that maps $H\times H$ to 1.
Then $\Hom(H,K)=T^{-1}(1)$ where the map
$$\displaylines{
  T:C(H,K)\to C(H\times H,K) \cr
  \noalign{\line{is defined as\hfill}}
  T(\sigma)(h_1,h_2)=\sigma(h_1)\,\sigma(h_2)\,\sigma(h_1h_2)^{-1}. \cr
               }
$$
Let $\sigma\in\Hom(H,K)$.
Then the $K$-orbit through $\sigma$ is given by the range
of the map
$$\displaylines{
  S:K\to C(H,K) \cr
  \noalign{\line{is defined as\hfill}}
  S(k)(h)=k\,\sigma(h)\,k^{-1}. \cr
               }
$$
The tangent space at $\sigma$ to the $K$-orbit through $\sigma$
is given by the range of the differential $(S_*)_1$.

We identify the tangent space at any point in the Lie group $K$
with the Lie algebra $\eufk$ by right translation.
If $\eufX$ is a compact topological space,
then this gives an identification of the tangent space
at any point in the Banach manifold $C(\eufX,K)$ 
with the Banach space $C(\eufX,\eufk)$.
Under these identifications, the differentials
$$\eqalignno{
    (S_*)_1&:\eufk\to C(H,\eufk) \cr
    (T_*)_\sigma&:C(H,\eufk)\to C(H\times H,\eufk) \cr
\noalign{\line{are given by\hfil}}
    (S_*)_1(\xi)(h)&=\xi-\Ad\sigma(h)\,\xi \cr
    (T_*)_\sigma(\lambda)(h_1,h_2)&=\lambda(h_1)+\Ad\sigma(h_1)\,\lambda(h_2)
        -\lambda(h_1h_2) . \cr
               }
$$
We see that $(S_*)_1=-\delta_1$ and $(T_*)_\sigma=\delta_2$ where
the linear maps
$$\delta_q:C(H^{q-1},\eufk) \to C(H^q,\eufk)$$
are defined as
$$\eqalign{
  \delta_q\mu(h_1,\ldots,h_q) & = \Ad\sigma(h_1)\,\mu(h_2,\ldots,h_q) \cr
    & \qquad  + \sum_{i=1}^{q-1} (-1)^i \mu(h_1,\ldots,h_ih_{i+1},\ldots,h_q)
              + (-1)^q \mu(h_1,\ldots,h_{q-1}) . \cr
          }
$$
A short calculation shows that the linear maps $\delta_q$ define
a chain complex.
The cohomology groups $H^q_\cont(H,\eufk)$
of this complex are called the {\it continuous cohomology} groups
of $H$ with coefficients in the $H$-module $\eufk$.
Note that continuous cohomology is defined the same way as
group cohomology,
except that one only considers continuous cochains.

For a compact topological group $H$
the continuous cohomology groups are essentially trivial.
In fact, define maps
$$\displaylines{
    s_q: C(H^q,\eufk) \to C(H^{q-1},\eufk) \cr
   \noalign{\line{by\hfill}}
    s_q \mu (h_1,\ldots,h_{q-1})
       = (-1)^q \int_H \mu (h_1,\ldots,h_{q-1},h)\,dh , \cr
               }
$$
where $dh$ is the normalized Haar measure on $H$.
Then a short calculation shows that
$$\cases{
  s_1\delta_1 = 1-p \cr
  \delta_q s_q+s_{q+1}\delta_{q+1} = 1 & for $q\ge1$, \cr
        }
$$
where $p(\xi)=\int_H \Ad\sigma(h)\,\xi\,dh$.
The map $p$ is a projection of
$\eufk$ onto the fixed point set $\eufk^H$.
It follows that
$$H^q_\cont(H,\eufk) = \cases{ \eufk^H & for $q=0$ \cr 0 & for $q\ge1$ .}$$
We conclude that the null space of $(T_*)_\sigma$
is the tangent space of the $K$-orbit through $\sigma$
and that $(T_*)_\sigma$ has closed range.
It then follows from the implicit function theorem that
the $K$-orbit through $\sigma$ has a neighborhood in $C(H,K)$
where $T\ne1$ away from the $K$-orbit itself.
\enddemo

\proclaim Lemma 5.2.

If $H$ is a semisimple compact Lie group
and $K$ is a compact Lie group,
then there exist only finitely many $K$-conjugacy classes
of homomorphisms $H\to K$.

\demo Proof.

Let $\eufh$ and $\eufk$ be the Lie algebras of $H$ and $K$.
Let $\hom(\eufh,\eufk)$ be the space of linear maps
$\eufh\to\eufk$ and $\eufhom(\eufh,\eufk)$ the
set of Lie algebra homomorphisms $\eufh\to\eufk$.
Since $\eufh$ is semisimple,
any homomorphism $\eufh\to\eufk$ takes values in the semisimple
part $\eufk_\ss$ of $\eufk$.
The Killing form $|X|^2=-\Tr(\ad X)^2$
defines norms on $\eufh$ and $\eufk_\ss$.
This gives a norm
$$|\lambda|^2  = \sup_{0\ne X\in\eufh} {|\lambda X|^2  \over |X|^2 }
               = \sup_{0\ne X\in\eufh}
                    { \Tr(\ad(\lambda X))^2 \over \Tr (\ad X)^2}
$$
on $\hom(\eufh,\eufk_\ss)$.
Now $\ad\circ\lambda$ gives a representation of $\eufh$ on $\eufk$.
We see that $|\lambda|$ only depends on
(the isomorphism class of) this representation.
It follows from the classification of 
representations of semisimple Lie algebras
that there are only finitely many representations of $\eufh$
of given dimension.
Hence $|\lambda|$ can assume only finitely
many values.
In particular, $\eufhom(\eufh,\eufk)$ is a compact subset
of $\hom(\eufh,\eufk_\ss)$,
and hence a compact subset of $\hom(\eufh,\eufk)$.

If $\sigma\in\Hom(H,K)$,
then $(\sigma_*)_1\in\eufhom(\eufh,\eufk)$.
Thus we have a uniform bound for $(\sigma_*)_1$,
and hence for $(\sigma_*)_h$ for any $h\in H$.
It then follows from the Arzela-Ascoli theorem that
$\Hom(H,K)$ is a compact subset of $C(H,K)$.
Lemma 5.2 now follows from Lemma 5.1.
\enddemo

Next we consider a short exact sequence
$$1\longrightarrow L \longrightarrow K \mapright{\tau} H \longrightarrow 1
  \eqno(5.1)
$$
of {\it compact} Lie groups.
A smooth homomorphism $\sigma:H\to K$ which is a right inverse of $\tau$
is called a splitting of the sequence.
Two splittings are said to be equivalent
if they differ by conjugation by an element of $L$.

\proclaim Lemma 5.3.

Any $K$-conjugacy class in $\Hom(H,K)$
contains only finitely many equivalence classes
of splittings of the short exact sequence {\rm(5.1)}.

\demo Proof.

If $\sigma\in\Hom(H,K)$ is a splitting,
then the other splittings in the $K$-conjugacy class of $\sigma$
are of the form $k\sigma k^{-1}$
with $k\in\tau^{-1}(Z(H))$.
The splittings $\sigma$ and  $k\sigma k^{-1}$
are equivalent, i.e.~they are $L$-conjugates,
if and only if $k\in L\,Z(K)$.
It follows that there is a 1--1 correspondence
between the $L$-conjugacy classes of splittings
in the $K$-conjugacy class of $\sigma$
and the elements of the group $\tau^{-1}(Z(H))/L\,Z(K)\cong Z(H)/\tau(Z(K))$.
The Lie algebras of $Z(H)$, $\tau(Z(K))$ and
$Z(H)/\tau(Z(K))$ are
$\eufz(\eufh)$, $\tau_*(\eufz(\eufk))$ and
$\eufz(\eufh)/\tau_*(\eufz(\eufk))$.
It follows from the structure theorem
for compact Lie groups that $\eufh=\eufz(\eufh)\oplus\eufh_\ss$
and $\eufk=\eufz(\eufk)\oplus\eufk_\ss$.
We have $\tau_*(\eufk_\ss)\subseteq\eufh_\ss$.
As $\tau_*$ is surjective,
we have $\tau_*(\eufz(\eufk))\subseteq\eufz(\eufh)$.
Hence
$\tau_*(\eufk_\ss)=\eufh_\ss$
and
$\tau_*(\eufz(\eufk))=\eufz(\eufh)$.
Hence $Z(H)/\tau(Z(K))$ is a discrete group.
It is a compact group and is hence a finite group.
The Lemma follows.
\enddemo

The following Lemmas are immediate consequences of Lemmas 5.1, 5.2 and 5.3.
The Lie group $L$ acts on $C(H,K)$ by conjugation.

\proclaim Lemma 5.4.

The splittings of the short exact sequence {\rm(5.1)}
form a closed discrete union of $L$-orbits in $C(H,K)$.

\proclaim Lemma 5.5.

If $H$ is semisimple,
then there exist only finitely many
equivalence classes of splittings of the short exact sequence {\rm(5.1)}.


\heading \S6. The fibers of the map $\pi:\eusE\to\eusB^H$

The main ingredient in the proof of the reverse
implication in Theorem 2 is the following Proposition:

\proclaim Proposition 6.1.

For any $[A]\in\eusB^H$,
$\pi^{-1}[A]$ is a closed discrete subset of $\eusE$.

\demo Proof.

The set is obviously closed. We will now show that it is discrete.
Any sequence in $\pi^{-1}[A]$ is of the form
$[g_n.A,\sigma]$ with $g_n\in\eusG$ and
$$g_n.A\in\eusA^H .
  \eqno(6.1)
$$
Assume that that $[g_n.A,\sigma]\to[A_\infty,\sigma]$
in $L^{1,2}_3$ for some $A_\infty\in\eusA^H$.
This means that there exists a sequence
$$\gamma_n\in\eusG^H
  \eqno(6.2)
$$
such that
$$\gamma_n g_n.A\to A_\infty
  \rlap{\qquad in $L^{1,2}_3$\quad as $n\to\infty$.}
  \eqno(6.3)
$$
We have to show that $[g_n.A,\sigma]=[A_\infty,\sigma]$ for large $n$.

The orbits of the action of $\eusG$ on $\eusA$ are closed,
so $A_\infty$ lies in the $\eusG$-orbit of $A$.
Then there exists a sequence $k_n\in\eusG$ such that
$$k_n\to1 \rlap{\qquad in $L^{2,2}_3$\quad as $n\to\infty$}
  \eqno(6.4)
$$
and
$k_n \gamma_n g_n.A=A_\infty$.
It then follows from the compactness of the isotropy subgroups
that, after passing to a subsequence,
$k_n \gamma_n  g_n\to g_\infty$ in $L^{2,2}_3$ as $n\to\infty$
for some $g_\infty\in\eusG$.
It then follows from (6.4) that
$$\gamma_n  g_n\to g_\infty
    \rlap{\qquad in $L^{2,2}_3$\quad as $n\to\infty$}
  \eqno (6.5)
$$
and from (6.3) that
$$g_\infty.A=A_\infty .
  \eqno(6.6)
$$
By (6.1),
$g_ng_\infty^{-1}.A_\infty=g_n.A\in\eusA^H$.
In other words, $g_ng_\infty^{-1}.A_\infty$ is invariant under $\sigma$,
so $A_\infty$ is invariant under $g_\infty g_n^{-1}.\sigma$.
A lift that preserves $A_\infty$ is by definition the same as
a continuous right inverse of the homomorphism $\StabB{A_\infty}\to H$.
Thus $g_\infty g_n^{-1}.\sigma$ is a homomorphism
$H\to\StabB{A_\infty}$,
in particular an element of $C(H,\StabB{A_\infty})$.
By (6.2) and (4.1),
$\gamma_n.\sigma=\sigma$.
It then follows from (6.5) that
$$g_\infty g_n^{-1}.\sigma
 = g_\infty (\gamma_n g_n)^{-1} .\sigma
 \to \sigma
  \hbox{\quad in $C(H,\StabB{A_\infty})$\quad as $n\to\infty$.}
$$
It follows from Lemma 5.4,
applied to the short exact sequence (4.4) with $A=A_\infty$,
that $g_\infty g_n^{-1}.\sigma$ lies in the same 
$\StabA{A_\infty}$-orbit as $\sigma$ for large $n$.
In other words,
there exists a sequence
$$s_n\in\StabA{A_\infty}
  \eqno(6.7)
$$
such that $s_n g_\infty g_n^{-1}.\sigma=\sigma$ for large $n$.
By (4.1),
$$s_n g_\infty g_n^{-1}\in\eusG^H
  \eqno(6.8)
$$
for large $n$.
It follows from (6.8), (6.6) and (6.7) that
$$[g_n.A,\sigma]
  = [s_n g_\infty.A,\sigma]
  = [s_n.A_\infty,\sigma]
  = [A_\infty,\sigma]
$$
for large $n$.
\enddemo

\demo Proof of Theorem 2, (1)$\Rightarrow$(2).

Assume that Condition (2) were false.
Then there would exists a Yang-Mills connection $[A]\in\eusB^H$ 
with infinitely many distinct preimages $[A_n,\sigma]$ in $\eusE$.
Then $A_n$ would be a Palais-Smale sequence.
It follows from Proposition~6.1 that the sequence $[A_n,\sigma]$ would
not contain any convergent subsequence.
In other words, Condition (1) would be false.
\enddemo

In order to state the next Proposition
we need to consider all lifts $\sigma$ of $\rho$ simultaneously.
Thus we choose one representative $\sigma_i$ from each
gauge equivalence class of lifts of $\rho$;
here $i$ ranges over some set $I$
that parametrizes the gauge equivalence classes of lifts.
Each lift $\sigma_i$ gives actions
$H\times\eusA\to\eusA$ and $H\times\eusG\to\eusG$.
We let $\eusA^H_i$ and $\eusG^H_i$ be the fixed point sets
of these actions and we let $\eusE_i=\eusA^H_i/\eusG^H_i$.
We denote the equivalence class of $A\in\eusA^H_i$ by
$[A,\sigma_i]\in\eusE_i$.
There are natural maps 
$$\pi_i:\eusE_i\to\eusB^H.$$

\proclaim Proposition 6.2.

Let $A\in\eusA$ with $[A]\in\eusB^H$.
Then there is a 1--1 correspondence between the disjoint union
$\coprod_{i\in I}\pi_i^{-1}[A]$ and the set of equivalence classes 
of splittings of the short exact sequence {\rm(4.4)}.

\demo Proof.

It follows from (4.1) that there is a 1--1 correspondence
$$\coprod_{i\in I}\eusE_i \leftrightarrow
  \bigg\{\> (A,\sigma') \>\bigg|\>  \vcenter{
       \hbox{$A\in\eusA$, $\sigma'$ is a lift of $\rho$,}
       \hbox{and $A$ is $\sigma'$-invariant}
                                        } \>\bigg\}
  \bigg/ \eusG .
  \eqno(6.9)
$$
We denote the equivalence class of $(A,\sigma')$ by $[A,\sigma']$.
That is consistent with the notation $[A,\sigma_i]$
for the equivalence class of $A\in\eusA^H_i$ in $\eusE_i$.
In other words, the map to the right in (6.9) is given by
$[A,\sigma_i]\mapsto[A,\sigma_i]$.

Let $A\in\eusA$ with $[A]\in\eusB^H$.
Under the 1--1 correspondence (6.9)
the preimage of $[A]$ in $\coprod_{i\in I}\eusE_i$
corresponds to the set of all equivalence classes $[A,\sigma']$
such that $A$ is invariant under $\sigma'$.
A lift $\sigma'$ that preserves $A$ is by definition
the same as a splitting of the 
the short exact sequence (4.4).
We have $[A,\sigma']=[A,\sigma'']$
if and only if there exists $g\in\eusG$
such that $g.A=A$, i.e.~$g\in\StabA{A}$, and $g.\sigma'=\sigma''$.
In other words, two splittings $\sigma'$ and $\sigma''$
of (4.4) give the same element
of $\coprod_{i\in I}\eusE_i$ if and only if they are $\StabA{A}$-conjugates.
\enddemo

Theorem 3(a) follows from Proposition 6.2.
Theorem 3(b) follows from Proposition~6.2 and Lemma 5.5.


\references

\ref[1]
D.R. Adams,
{\it A note on Riesz potentials},
Duke Math. J. {\bf42} (1975) 765--778.

\ref[2]
G. Bor,
{\it Yang-Mills fields which are not self-dual},
Comm. Math. Phys. {\bf145} (1992) 393--410.

\ref[3]
P.J. Braam and G. Mati\'c,
{\it The Smith conjecture in dimension four and equivariant gauge theory},
Forum Math. {\bf5} (1993) 299--311.

\ref[4]
S. Campanato,
{\it Propriet\`a di inclusione per spazi di Morrey},
Ricerche Mat. {\bf12} (1963) 67--86.

\ref[5]
\vrule height 0pt depth 0.3pt width 20pt \thinspace ,
{\it Equazioni ellittiche del $\hbox{\it II}^{\,\circ}$ ordine e spazi
$\eusL^{(2,\lambda)}$},
Ann. Mat. Pura Appl. (4) {\bf69} (1965) 321--381.

\ref[6]
F. Chiarenza and M. Frasca,
{\it Morrey spaces and Hardy-Littlewood maximal function},
Rend. Mat. Appl. (7) {\bf7} (1987) 273--279.

\ref[7]
W.T. van Est,
{\it On the algebraic cohomology concepts in Lie groups I and II},
Indag. Math. {\bf17} (1955) 225--233, 286--294.

\ref[8]
M. Furuta,
{\it A remark on a fixed point of finite group action on $S^4$},
Topology {\bf28} (1989) 35--38.

\ref[9]
M. Giaquinta,
{\it Multiple integrals in the calculus of variations and nonlinear
elliptic systems},
Princeton Univ. Press, Princeton, 1983.

\ref[10]
D. Gilbarg and N.S. Trudinger,
{\it Elliptic partial differential equations of second order},
$2^{\rm nd}$ ed.,
Springer-Verlag, Berlin, 1983.

\ref[11]
U. Gritsch,
{\it Morse theory for the Yang-Mills functional over equivariant
four-manifolds and equivariant homotopy theory},
Ph.D. thesis, Stanford University, 1997.

\ref[12]
\vrule height 0pt depth 0.3pt width 20pt \thinspace ,
{\it Morse theory for the Yang-Mills functional 
via equivariant homotopy theory},
preprint, University of Bonn, 1997.

\ref[13]
U. Gritsch and J. R\aa de,
{\it Equivariant Morse theory for invariant $\SU(2)$-connections 
over $S^2\times S^2$},
in preparation.

\ref[14]
A. Guichardet,
{\it Cohomologie des groupes topologiques et des alg\`ebres de Lie},\break
Cedic/Nathan, Paris, 1980.

\ref[15]
R. Haas,
Ph.D. thesis, University of Texas at Austin, in preparation.

\ref[16]
G. Hochschild and G.D. Mostow,
{\it Cohomology of Lie groups},
Illinois J. Math. {\bf6} (1962) 367--401.

\ref[17]
D.H. Lee,
{\it On deformations of homomorphisms of locally compact groups},
Trans. Amer. Math. Soc. {\bf 191} (1974) 353--361.

\ref[18]
C.B. Morrey, Jr.,
{\it Multiple integral problems in the calculus of variations
and related topics},
Univ. California Publ. Math. (N.S.) {\bf1} (1943) 1--130.

\ref[19]
G.D. Mostow,
{\it Cohomology of topological groups and solvmanifolds},
Ann. of Math. (2) {\bf73} (1961) 20--48.

\ref[20]
A. Nijenhuis and R.W. Richardson, Jr.,
{\it Deformations of homomorphisms of Lie groups and Lie algebras},
Bull. Amer. Math. Soc. {\bf73} (1967) 175--179.

\ref[21]
T.H. Parker,
{\it A Morse theory for equivariant Yang-Mills},
Duke Math. J. {\bf66} (1992) 337--356.

\ref [22]
J. Peetre,
{\it On convolution operators leaving $L^{p,\lambda}$ spaces invariant},
Ann. Mat. Pura Appl. (4) {\bf 72} (1966) 295--304.

\ref[23]
L. Sadun and J. Segert,
{\it Non-self-dual Yang-Mills connections with quadrupole symmetry},
Comm. Math. Phys. {\bf145} (1992) 363--391.

\ref[24]
L.M.~Sibner, R.J.~Sibner and K.~Uhlenbeck,
{\it Solutions to Yang-Mills equations that are not self-dual},
Proc. Nat. Acad. Sci. U.S.A. {\bf86} (1989) 8610--8613.

\ref [25]
M.E. Taylor,
{\it Analysis on Morrey spaces and applications to Navier-Stokes and
other evolution equations},
Comm. Partial Differential Equations {\bf17} (1992) 1407--1456.

\ref [26]
H. Urakawa,
{\it Equivariant theory of Yang-Mills connections over Riemannian
manifolds of cohomogeneity one},
Indiana Univ. Math. J. {\bf37} (1988) 753--788.

\ref[27]
H.-Y. Wang,
{\it The existence of nonminimal solutions to the Yang-Mills equation
with group $\SU(2)$ on $S^2\times S^2$ and $S^1\times S^3$},
J. Differential Geom. {\bf34} (1991) 701--767.

\ref[28]
A. Weil,
{\it Remarks on the cohomology of groups},
Ann. of Math. (2) {\bf80} (1964) 149--157.

\endlines{Department of Mathematics, G\"oteborg University,
  S-412\thinspace96 G\"oteborg, Sweden\cr
  http://www.math.chalmers.se/$\sim$jrade/\cr
  jrade@math.chalmers.se\cr
  August 1, 1997\cr
         }
 
\bye